\documentclass[11pt]{amsart}
\usepackage{amssymb,amscd,epsfig}

\def\e{{\varepsilon}}
\def\d{{\delta}}
\def\bd{{\rm bd}}

\begin{document}

\newtheorem{defi}{Definition}
\newtheorem{theorem}{Theorem}
\newtheorem{remark}{Remark}
\newtheorem{lemma}{Lemma}
\newtheorem{coro}{Corollary}
\newtheorem{quest}{Question}

\title{Total curvature and spiralling  shortest paths.}

\author{I.~B\'{a}r\'{a}ny}
\address{R\'enyi Institute of Mathematics,
POB 127, 1364 Budapest,
Hungary, and
Department of Mathematics,
University College London,
Gower Street London, 
WC1E6BT United Kingdom}
\email{barany@renyi.hu}

\author{ K.~Kuperberg}
\address{Department of Mathematics, Auburn University, Auburn, AL 36830-5310, USA}
\email{kuperkm@math.auburn.edu}

\author{ T.~Zamfirescu}
\address{ Fachbereich Mathematik,
Universit\"at Dortmund,
44221 Dortmund, Germany}
\email{tudor.zamfirescu@mathematik.uni-dortmund.de}

\thanks{I.~B\'{a}r\'{a}ny's research was supported in part by Hungarian
National Science Foundation grants \#T-032452 and \#T-029255}
\thanks{K.~Kuperberg's research was supported in part by US National
Science Foundation grant \#DMS-9704558}
\begin{abstract} This paper gives a partial confirmation of a conjecture of 
P.~Agarwal, S.~Har-Peled, M.~Sharir, and K.~Varadarajan that the total curvature of a
shortest path on the boundary of a convex polyhedron in ${\mathbb R}^3$ cannot be
arbitrarily large.  It is shown here that the conjecture holds for a class of polytopes
for which the ratio of the radii of the circumscribed and inscribed ball is bounded. On
the other hand, an  example is constructed to show that the total curvature of a 
shortest path on the boundary of a convex polyhedron in ${\mathbb R}^3$ can
exceed $2\pi$. Another example shows that the spiraling number of a 
shortest path on the boundary of a convex polyhedron can be arbitrarily
large.
\end{abstract}
\maketitle

\section{Introduction} The {\em total curvature\/} of a $C^2$ path $C$
parameterized by  arc length $s$ in ${\mathbb R}^n$, $r(s)$,  is defined as
$\int _C |r''(s)|\,  ds$. W.~Fenchel proved in 1929 for ${\mathbb R}^3$ and
K.~Borsuk in 1947 for any ${\mathbb R}^n$ that  the total curvature of a
closed curve is bounded from below by $2\pi$, with the equality holding only
for convex simple closed curves in ${\mathbb R}^2$. The {\em total
curvature\/} of a polygonal path $P=[z_0,z_1,\ldots , z_n]$ is defined as 
$$t(P)=\sum_{i=1}^{n-1}( \pi -\angle z_{i-1} z_i z_{i+1}) .$$ 

Let $\mathcal K$ be the set of all compact convex polyhedra in ${\mathbb
R}^3$. Let $\mathcal T=\{t(P)\}$, where $P$ is a shortest path joining two
points on the boundary of a polyhedron  $K\in {\mathcal K}$.  It has been
asked in \cite{Agarwal}  whether the set $\mathcal T$ is bounded. 

We shall prove here, in Theorem~\ref{rB} below, that the conjecture holds
for polytopes $K$ for which the ratio $R/r$ is bounded from above: here $R$
and $r$, respectively, are the radii of the circumscribed and inscribed ball
to $K$.

We define the  {\em spiralling number} $s(P)$ for the path $P$  from $a$ to
$b$ on the polytopal surface by  considering a variable point $x\in
P\setminus\{a,b\}$, writing this point in cylindrical coordinates as
$(2\pi\phi(x), r(x), z(x))$  where
the $z$-axis is the line through $a$ and $b$, and $\phi$ is a
continuous function. Set now  
$$s(P)=\lim_{x\to a, y\to b}|\phi(x)-\phi(y)|.$$ This measures how many
times the path spirals around the line through $a$  and $b$.

The proof method of Theorem~\ref{rB} would work for all polytopes if the
function $\phi$ had bounded variation. This, however, will be shown to be
false. In section~\ref{spiral}, we will construct needle-like
polytopes $K$ with shortest
path $P$ such that $s(P)$ is arbitrarily large. Even more surprisingly, the
example can be modified so that $P$ spirals around the line through $a$ and
$b$ a hundred times in one direction, then two hundred times in the opposite
direction, then thousand times in the first direction, etc.

The total curvature $t(P)$ of a planar path $P$ is bounded sharply by $2\pi$ and
this bound is the lowest possible.   A  triangle with one of the angles very
close to $\pi$ and two points on the two sides adjacent to the wide angle
but close to the vertices at the acute angles provides a simple example.

In section~\ref{total}, we construct an example of a shortest path on the
boundary of a convex polyhedron in ${\mathbb R}^3$ for  which the $2\pi$
bound does not hold. On some  level the example resembles the  planar example 
involving a triangle with one of the angles close to $\pi$.

\section{Bounded total curvature of shortest paths for $R/r$ bounded}

On the boundary of a polytope consider 
a shortest path  $P$ with non-smooth points $z_0,z_1,...,z_n$.
Put $x_i=(z_i-z_{i-1})/\|z_i-z_{i-1}\|$\ 
$(i=1,...,n)$. Let $u_i$ be the outernormal 
(of the interior) of the facet of $z_{i-1}$ and $z_i$, and let $\xi_i$ be the angle
between $x_i$ and $x_{i+1}$\ $(i=1,...,n-1)$.
Then the total curvature of $P$ is 
$\sum_{i=1}^{n-1}\xi_i$.
This can be easily checked.
We remark that $x_i-x_{i+1}=\lambda_i(u_i + u_{i+1})$ with $\lambda_i>0.$ 

\begin{lemma}\label{unit}  If there is a unit vector v such that
$u_iv\ge\eta>0$ for  all i, then $\sum_{i=1}^{n-1}\xi_i<\frac{\pi}{\eta}.$
\end{lemma}

\begin{proof}
First remark that $\xi<\frac {\pi}{2}\sin \xi$ if $\xi<\pi/2$. Then
$$\sum_{i=1}^{n-1}\xi_i= 2\sum_{i=1}^{n-1}\frac{\xi_i}{2}< 2\sum_{i=1}^{n-1}
\frac{\pi}{2} \sin \frac{\xi_i}{2} $$
$$=\frac{\pi}{2}\sum_{i=1}^{n-1}\|x_i-x_{i+1}\|=\frac{\pi}{2}\sum_{i=1}^{n-1}
\|\lambda_i(u_i+u_{i+1})\|.$$
Since $\|u_i+u_{i+1}\|\le 2$ and $(u_i+u_{i+1})v\ge 2\eta$, we have
$$\|u_i+u_{i+1}\|\le \frac{(u_i+u_{i+1})v}{\eta}.$$
Hence 
$$\sum_{i=1}^{n-1}\xi_i<\frac{\pi}{2}
\sum_{i=1}^{n-1}\lambda_i\frac{(u_i+u_{i+1})v}{\eta}
=\frac{\pi}{2\eta}\sum_{i=1}^{n-1}(x_i-x_{i+1})v
=\frac{\pi}{2\eta}(x_1-x_n)v\le\frac{\pi}{\eta}.$$
\end{proof}

Denote by $B$ the closed unit ball in ${\mathbb R}^3$.

\begin{theorem}\label{rB} 
Let a polytope $Q$ satisfy $rB\subset Q\subset B$. Then the total curvature
of any shortest path on the boundary of Q is less than $4\pi^2r^{-2}.$
\end{theorem}

\begin{proof}
Let $P$ be a shortest path between points $a$ and $b$ on the boundary  bd$Q$
of  $Q$. The length of $P$ is less than $\pi$. Indeed any plane through $a$
and $b$  intersects the boundary of $B$ along a circle of length at most
$2\pi$, and  the boundary of $Q$ along a polygon of even smaller length. So
one of the  two broken lines into which $a$ and $b$ divide the polygon must
have length less than $\pi$, and $P$ cannot be longer.  For a non-zero vector
$v\in {\mathbb R}^3$ define the $v$-shadow $S(v)$ of $Q$ as
$$S(v)= \bd Q\cap\{\frac{r}{2}B + \lambda v\ |\ \lambda\ge 0\}.$$ Assume that $u$ is
the outernormal at an interior point $f$ of a facet $F$ of
$Q$ with $f\in S(v)$. We claim that $ {uv}/{\|v\|}>r/2.$ Indeed, the plane
$\Pi$ of $F$ does not meet  the interior of $rB$. Let $\gamma$ be the angle
between $u$ and $v$. Of course, $\|f\|<1$. The distance from $f$ to the line
through the origin  $\bf 0$ and $v$ is at most $r/2$ because $f \in
S(v)$. Further, the distance from
$\bf 0$ to $\Pi$ is at least $r$. Therefore
$$\frac{ uv}{\|v\|} = \cos \gamma > \frac{r}{2},$$  
which proves the claim. We are going to define points $v_1, v_2,..., v_k$ on
$P$ by recursion. Set  $v_1=a$ and assume that $v_i$ has been defined. Let
$v_{i+1}$ be the first point  on $P$, going from $v_i$ to $b$, which is not
in the interior of $S(v_i)$. Write $P_i$ for the part of the path $P$
between $v_i$ and $v_{i+1}$. Since the length of $P_i$ is at least
$\|v_i-v_{i+1}\|\ge r/2$ and the length of $P$ is less than
$\pi$, there are  only  $k< 2\pi/r$ paths $P_i$, and we put $v_{k+1}=b$. 
Note that Lemma~\ref{unit} applies to $P_i$ with $\eta=r/2$ by the above
claim. Hence the total curvature of $P_i$ is less than $2\pi/r$. Summing up,
the total  curvature of $P$ is less than
$$k\frac{2\pi}{r}< \frac{4\pi^2}{r^2}.$$
\end{proof}

Now we formulate Theorem~\ref{rB} differently. Let $E$ be the ellipsoid  of
largest volume included in the polytope $Q$. As the total curvature is
invariant under angle preserving linear transformation, we may assume that
$E$  has half-axes $a,b,c$ with $0<a\le b\le c=1.$ Call $Q$ {\it
needle-like} if $b$ is small, and {\it pancake-like} if $a$ is small
compared with $b$. Theorem~\ref{rB} on the 
boundedness of total curvature holds for convex bodies that are not 
pancake-like. Only moderate effort is needed to prove the following: If the
conjecture holds for needle-like convex bodies, then it holds for
pancake-like ones as well. So it would be enough to prove the conjecture for
needle-like convex bodies. This would follow if the total variation of
$\phi$, the function defined in  the first paragraph, were bounded. But this
is not true:  an example showing this is in section~\ref{spiral}.

\section{The total curvature $t(P)>2\pi$: an example}\label{total} 

We construct a convex body $\Delta$ with two points on the
boundary such that the total curvature of the shortest path joining 
the points exceeds $2\pi$. $\Delta$ is constructed in four steps, with
certain unbounded convex bodies $U$, $X$, and $Y_\alpha$ described in
steps 1-3. The Cartesian coordinates of a point $x\in {\mathbb R}^n$ are
denoted by $x^{(1)}, \ldots , x^{(n)}$.

For $0<\epsilon < \frac{\sqrt{2}}{8\pi}$  and $i=1,2$,  define $V_i$ to be the
plane $\{ (x^{(1)},x^{(2)},x^{(3)})\in {\mathbb R}^3\ |\
x^{(3)}=(-1)^i\epsilon\}$. Let $A_i$ be the parabola in 
$ \{(x^{(1)},x^{(2)},x^{(3)})\in V_i\ |\  x^{(2)}=(x^{(1)})^2\}$  
and $B_i$ the parabola in 
$\{ (x^{(1)},x^{(2)},x^{(3)})\in V_i\ |\  x^{(2)}=2  (x^{(1)})^2\}$. 
Denote by $v_i$ the common vertex of $A_i$ and  $B_i$, by   $a_i$ the focus
of $A_i$, and  by  $b_i$ the focus of $B_i$. Thus 
$a_i=(0, \frac{1 }{4},(-1)^i\epsilon )$, 
$b_i=(0, \frac{1}{8} ,(-1)^i\epsilon )$, and  $v_i=(0,0,(-1)^i\epsilon )$.
Denote by $R_i$ the convex region in the plane $V_i$ bounded by $A_i$.  Put
$$L=\{ (x^{(1)},x^{(2)},x^{(3)})\in {\mathbb R}^3\ |\
(x^{(2)}-\epsilon)^2+(x^{(3)})^2=2\epsilon^2,\ x^{(1)}=0,\ 
x^{(2)}\leq 0\}$$ and
$$\Gamma=\{ (x^{(1)},x^{(2)},x^{(3)})\in {\mathbb R}^3\ |\
( x^{(2)}-\epsilon)^2+(x^{(3)})^2=2\epsilon^2,\
-1/2\leq x^{(1)}\leq 1/2,\ x^{(2)}\leq 0\}.$$
Thus $L$ is a quarter-circle joining the vertices 
$v_1$ and $v_2$ and  $\Gamma$ is a surface  obtained by sliding $L$ along the segment
$\{ (x^{(1)},x^{(2)},x^{(3)})\in {\mathbb R}^3\ |\
 -1/2\leq x^{(1)}\leq 1/2,\ x^{(2)}=x^{(3)}=0 \}$. Note that
$|L|=\frac{\pi\epsilon}{\sqrt{2}}< \frac{1}{8}$.
Let $$U={\rm conv}(A_1\cup A_2\cup \Gamma)$$ be the  convex hull of the union
$A_1\cup A_2\cup \Gamma.$

Suppose that $G$ is a line in $V_2$, crossing $A_2$, and parallel to the
directrix of $A_2$.  Denote by $\delta_B$ the distance between $G$ and the
directrix of $B_2$, and by  $\delta_A$ the distance between $G$ and the
directrix of $A_2$. Note that if a shortest path $P$ in bd$U$ joining $G$
with $b_1$ passes through $v_1$ and $v_2$, then the length $|P|$ of $P$
equals $\delta_B+|L| $.

\begin{lemma}\label{path} The shortest path $P$ on the boundary 
of $U$ joining $G$ with
$b_1$ is unique and passes through the vertices $v_1$ and $v_2$.
\end{lemma}

\begin{proof} 
Suppose that $P'$ is a shortest path in {\rm bd}$U$ joining $b_1$
with $G$, different from $P$, crossing $A_1$ at
$d=(d^{(1)},d^{(2)},d^{(3)})$ and $A_2$ at $e=(e^{(1)},e^{(2)},e^{(3)})$,
see Fig.~\ref{paradist}. Note that  $d^{(2)}<e^{(2)}$.  

\begin{figure}[ht]
\centerline{\epsfig{figure=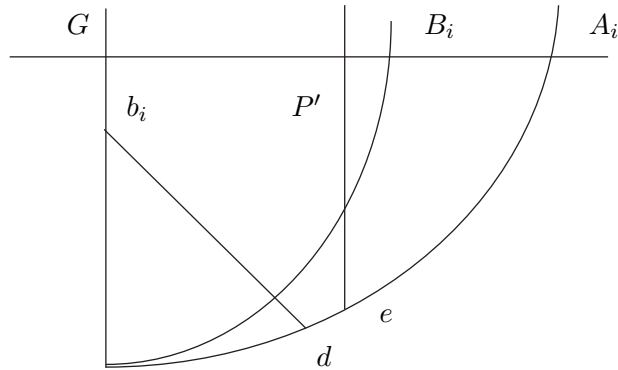,height=2.25in}}
\vspace{-2.05in}
\centerline{\hspace{0.25in}$G$\hspace{1.7in} $B_i$\hspace{0.65in} $A_i$}
\vspace{.25in}
\centerline{ $b_i$\hspace{.7in} $P'$\hspace{1in}$\ $}
\vspace{.9in}
\centerline{ $\ $\hspace{.6in}$e$}
\vspace{.05in}
\centerline{ $d$}
\vspace{.2in}
\caption{The parabolas.}
\label{paradist}
\end{figure}

To show that $|P'| >|P|$, consider the two cases: $|e^{(1)}|\geq \frac{1}{2}$
and $|e^{(1)}|< \frac{1}{2}$. We have

\begin{enumerate}

\item If  $|e^{(1)}|\geq \frac{1}{2}$, then $e^{(2)}\geq \frac{1}{4}$ and
$|P'|> {\rm dist}(G, e)  +{\rm dist}(e, a_2) =\delta_A=\delta_B
+\frac{1}{8}> \delta_B+|L| =|P|$.

\item If  $|e^{(1)}|< \frac{1}{2}$, then $P'$ projects onto  $P$ decreasing
its length.
\end{enumerate}
\end{proof}

\begin{remark}\label{monotone} If $P'(e)$ denotes the shortest path in {\rm
bd}$U$ joining $b_1$ with $G$ and passing through $e\in A_2$, then the length
$|P'(e)|$ is a monotone function  of $e^{(1)}$ for $0< e^{(1)}< \frac{1}{2}$
(and for $-\frac{1}{2}< e^{(1)}< 0$).
\end{remark}

In the next step of the construction, we modify $U$ to obtain a convex
unbounded slab $X$. Let $L'$ be a convex curve in the $(x^{(1)},x^{(2)})$-plane
close to the segment  $(x^{(1)},0,0)$, $-1/2\leq x^{(1)}\leq 1/2$, and
$\Gamma '$  a positively curved surface in $U$, close to
$\Gamma$, obtained by sliding $L$ along $L'$, keeping $L$ parallel to the
$(x^{(2)},x^{(3)})$-plane, and the midpoint of $L$ and $L'$. 
Denote by $X$ the convex hull  ${\rm conv}(A_1\cup
A_2\cup \Gamma ')$. We require of $L'$ (and consequently of $\Gamma '$) that

\begin{enumerate}
\item  $\Gamma '$ contains the quarter-circle $L$;

\item the path $P$ is the shortest path joining $G$ with
$b_1$ on the boundary of $X$.
\end{enumerate}

\begin{figure}[ht]
\centerline{\epsfig{figure=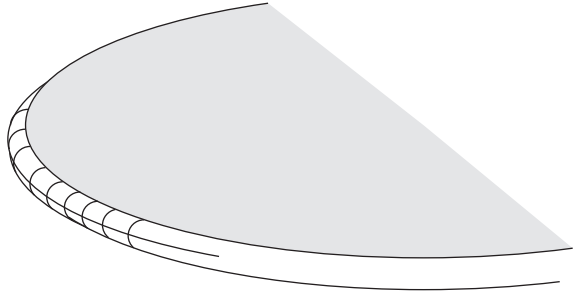,height=1.75in}}
\caption{The unbounded slab $X$.}
\vspace{.2in}
\label{slab}
\end{figure}

By Remark~\ref{monotone}, such an $X$ exists.  
Note that $X$ is a thin convex unbounded slab whose top and bottom are
(horizontal) planar  regions $R_i'$'s containing $R_i$'s and with
positively curved side surface close to the vertices $v_1$ and $v_2$,
see Fig.~\ref{slab}. 

For $0<\alpha <\frac{\pi}{4}$, let $Y_\alpha$ be a `slant slab'
obtained from $X$ by `$\alpha$-slanting' the plane $V_1$ at the
vertex $v_1$ in the direction perpendicular to the $x^{(1)}$-axis
as follows:  
Let $\Omega_\alpha (x)=(x^{(1)},x^{(2)}\cos \alpha , -\epsilon
-x^{(2)}\sin \alpha )$ and let $Y_\alpha={\rm conv}(X\cup
\Omega_\alpha (R_1'))$. Note that the portion of $Y_\alpha$ above
the $(x^{(1)},x^{(2)})$-plane is identical to that of $X$.

\begin{lemma}\label{tilted} The shortest path $\widetilde{P}$ joining
$G$ with $b_1'=\Omega_\alpha (b_1)$ on the boundary of $Y_\alpha$
passes through the vertices $v_1$ and $v_2$.
\end{lemma}

\begin{proof}  A path $P'$ joining $G$ with $b_1'$ in 
bd$Y_\alpha$ can be mapped (keeping the first coordinate and the
distance to the line $(t,0,-\epsilon)$ unchanged) onto the boundary
of  $X$ to a path joining $G$ with $b_1$. Such a projection does not
increase the length of the path. Hence, if 
$P'\not=\widetilde{P}$, $|P'| >\delta_B +\pi\epsilon =|P|=|\widetilde{P}|$.
\end{proof}

\begin{remark} For $\widetilde{P}$ as defined in Lemma~\ref{tilted}, the
total curvature of $\widetilde{P}$, $t(\widetilde{P})$,  equals
$\pi -\alpha$.
\end{remark}

\begin{figure}[ht]
\centerline{\epsfig{figure=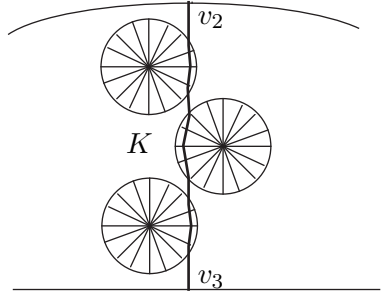,height=1.75in}}
\vspace{-1.7in}
\centerline{\hspace{0.25in}$v_2$}
\vspace{.5in}
\centerline{$K$\hspace{0.5in}}
\vspace{.5in}
\centerline{\hspace{0.25in}$v_3$}
\vspace{.2in}
\caption{Attached cones.}
\label{threecones}
\end{figure}

Denote by $S$ the side boundary of $Y_\alpha$, {\em i.e.}, $S={\rm bd}Y_\alpha
\setminus (R_2'\cup \Omega_\alpha (R_1'))$. For a point $w$, conv$(Y_\alpha ,w)$
denotes the convex hull of  $Y_\alpha$ and $w$. For a point
$w\notin Y_\alpha$,  define an {\em attached cone\/} with vertex $w$,
con$(w)$, as the closure of conv$(Y_\alpha ,w)- Y_\alpha$ provided that it 
does not intersect $R_2'\cup \Omega_\alpha (R_1')$.

Denote by $P_0$ the shortest path in $S$, a quarter-circle, joining  $v_2$ with
$v_3=(0, -\epsilon ,0)$. Let $w_1$, $w_2$, and $w_3$ be points close to
$P_0$ for which attached cones con$(w_i)$  are defined and are pairwise
disjoint. Note that $\widetilde{Y}=Y_\alpha \cup {\rm con}(w_1)
\cup {\rm con}(w_2)\cup {\rm con}(w_3)$ is convex. 

First choose points $w_1$, $w_2$, and $w_3$  so that

\begin{enumerate}
\item each con$(w_i)$, $i=1,2,3$, intersects $P_0$ at one point $t_i$ with 
${\rm dist} (t_1,t_2)={\rm dist} (t_2,t_3)$;

\item con$(w_1)$ and con$(w_3)$ are on the same side $P_0$, and con$(w_2)$
on the other side of $P_0$.
\end{enumerate}

Then move the points $w_1$, $w_2$, and $w_3$ slightly towards $P_0$ to
obtain points $w'_1$, $w'_2$, and $w'_3$ such that the attached cones
$C_i={\rm con}(w'_i)$ intersect $P_0$ and the shortest path $K$ joining $G$
and $b_1'$ in ${\rm bd}\widetilde{Y}$ crosses each $C_i$, but passes through
$v_1$,  $v_2$, and $v_3$, see Fig.~\ref{threecones}.

Note that $K$ is not planar and as it follows from the work of
W.~Fenchel~\cite{Fenchel} (alternatively, see K.~Borsuk~\cite{Borsuk} or
J.~Milnor~\cite{Milnor}) its total curvature is greater than that of $P_0$.
We have

\begin{lemma}\label{curvature}  The total curvature $t(K)= \pi +\beta
-\alpha$, where $\beta >0$.
\end{lemma}

\begin{figure}[ht]
\centerline{\epsfig{figure=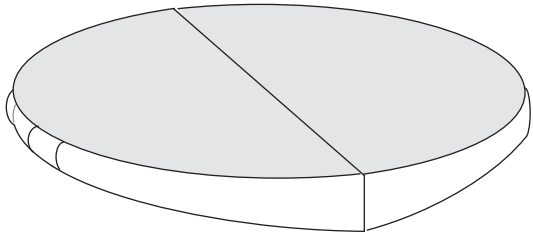,height=1.5in}}
\caption{The double slab $\Delta$.}
\label{doubleslab}
\vspace{.2in}
\end{figure}

Finally, take $\alpha < \beta $. The line $G$ consists of points $(t, g,
\epsilon )$, where $-\infty <t <\infty$. Let $\Delta _1$ be the part of
$\widetilde{Y}$ cut off by the plane   $x^{(2)}=g$. Let $\Delta _2$ be a
symmetrical copy of $\Delta _1$ and $\Delta =\Delta _1\cup \Delta _2$ with
$\Delta _1$ and $\Delta _2$ glued along the side $x^{(2)}=g$ (see
Fig.~\ref{doubleslab}). Let $K''$ and $b''_1$ be the path and the point in
$\Delta _2$ corresponding to $K'$ and $b_1'$. For sufficiently large $g$, the
path $\overline{K}= K'\cup K''$ is the shortest path  joining $b_1'$ and
$b''_1$ in $\Delta$. We have 
$T(\overline{K}) = 2(\pi -\beta + \alpha) > 2\pi$. Clearly $\Delta$ is
convex.

A polyhedral  example can be obtained by a suitable approximation of
$\Delta$.

\begin{theorem}\label{Delta} There exist  a convex polyhedron $M\subset
{\mathbb R}^3$ with two points $x$ and $y$ on the boundary  of $M$  such
that the total curvature of the (unique) shortest path joining $x$ and $y$
exceeds $2\pi$. 
\end{theorem}

\section{Spiralling shortest paths}\label{spiral}

Here we construct polytopal surfaces possessing shortest paths of
arbitrarily  large spiralling number. The intrinsic metric on these surfaces
will be denoted  by $\delta$.

\begin{theorem}\label{spiralling} Let $n$ be an integer. There exist a convex
polytope $Q$ and a shortest path $P$ between two points on  the boundary of
$Q$ with  $s(P)\ge n$.
\end{theorem}

\begin{proof}  Our example is the boundary of the convex hull of a family of
equilateral  triangles each two of which have pairwise parallel edges.
Suppose the construction performed up to the equilateral triangle $abc$,  so
that, from a  fixed point $x_0$ of the (already constructed) surface, the
intrinsic distance $\d(x_0, bc)$ to $bc$ (which is $\min_{y\in
bc}\d(x_0,y)$)  is smaller than the  distance to any other side of $abc$.
Let $\e$ be the smaller of the two differences. We now construct the next
triangle  $a'b'c'$. We do so that
$$\d(x_0,a'b')<\min\{\d(x_0,b'c'), \d(x_0, c'a')\}$$
and any shortest path from $x_0$ to any point $w$ of 
$a'b'$ necessarily crosses $bc$.
Let $z'$ be such that conv$\{z',a,b,c\}$ is a regular pyramid containing 
the already constructed surface. Then we can find a point $z$ behind $z'$
(see Fig.~\ref{pyramid}) but close to the pyramid axis such that 
$$\|z-a\|=\|z-b\|<\|z-c\|$$
and conv$\{z,a,b,c\}$ includes conv$\{z',a,b,c\}$. We shall choose $q>3$ and
put $a'=z+q(a-z)$, $b'=z+q(b-z)$ and $c'=z+q(c-z).$  The angles
$\alpha=\angle abb'$ and $\gamma=\angle cbb'$ satisfy
$$\pi/2<\gamma<\alpha$$
if $z$ is close enough to the pyramid axis. 

\begin{figure}[ht]
\centerline{\epsfig{figure=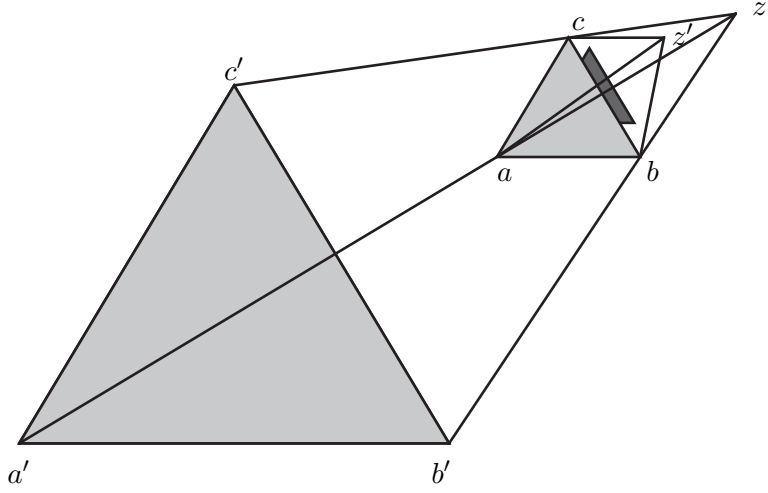,height=2.5in}}
\vspace{-2.55in}
\centerline{\hspace{4in}$z$}
\vspace{-.1in}
\centerline{\hspace{2.5in}$c$\hspace{.4in}}
\vspace{-.1in}
\centerline{\hspace{3.2in}$z'$}
\vspace{0in}
\centerline{\hspace{.5in}$c'$\hspace{2in}}
\vspace{.35in}
\centerline{\hspace{2.1in}$a$\hspace{.7in}$b$}
\vspace{1.4in}
\centerline{ $a'$\hspace{2.1in}$b'$\hspace{1.6in}}
\vspace{.2in}
\caption{Stacking triangles.}
\label{pyramid}
\vspace{.2in}
\end{figure}

Consider an unfolding of the
surface on the plane of $a, b, a', b'$ without cutting along $bb'$, but
cutting along $aa'$, and keep the notation. Let $s$ be the orthogonal
projection of $c$ on the line $L$ through $a$ and $b$. If $\alpha+\gamma-\pi$
is small enough, which is the case if $z$ is far  enough, then  $\d(c, s)$ 
is as small as desired, in particular smaller than $\e/2$. 
Let $v\in ab,  s'\in bs$, see Fig.~\ref{unfold}.

It is easily verified that $\delta(s',w) - \delta(v,w)$
increases when $w$ moves from $b'$ to $a'$. Thus
$$\delta(s',w) - \delta(v,w)\le \delta(s',a') - \delta(v,a').$$
If $\alpha+\gamma-\pi$ is small enough and $q$ large enough, then
$$\delta(s,a') - \delta(a,a')<\frac{\e}{2}.$$
Since $\delta(v,a')\ge\delta(a,a')$ and $\delta(s',a')\le\delta(s,a')$,
$$\delta(s',a') - \delta(v,a')<\frac{\e}{2}.$$
Hence
$$\delta(s',w)<\delta(v,w)+\frac{\e}{2}.$$
This ensures that any shortest path from $x_0$ to some point $w\in a'b'$ 
crosses  $bc$. Indeed, if the above path crosses  $ab$ at $v$, then
$$\d(x_0,v)+\d(v,w)\ge \d(x_0,ab)+\d(v,w)$$
$$\ge \d(x_0,bc)+\e+\d(v, w)
=\delta(x_0,u) + \frac{\e}{2} + \frac{\e}{2} + \delta(v,w)$$
$$>\d(x_0,u)+\d(c,s) +\delta(s',w)\ge \d(x_0,u)+\d(u,s')+\delta(s',w),$$
$u$ being a point of $bc$ closest to $x_0$, and  $s'$ the orthogonal
projection of $u$ on $L$; we got a contradiction. If the path crosses $ca$
at $v'$, say, then 
$$\d(x_0,v')+\d(v',w)\ge
\delta(x_0,ca)+\delta(v',w)\ge
\d(x_0, bc)+\e+\d(v,w)$$
and a contradiction is obtained as above.

\begin{figure}[ht]
\centerline{\epsfig{figure=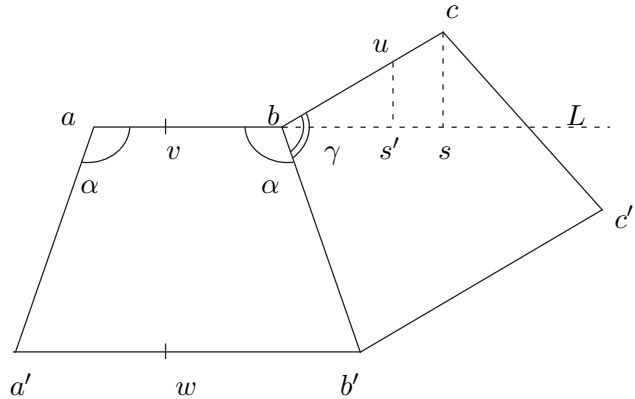,height=2in}}
\vspace{-2.1in}
\centerline{\hspace{1.6in}$c$}
\vspace{0in}
\centerline{\hspace{.85in}$u$}
\vspace{.2in}
\centerline{\hspace{.25in}$a$\hspace{1in}$b$\hspace{1.5in}$L$}
\vspace{0in}
\centerline{\hspace{.1in}$v$\hspace{.75in}$\gamma$\hspace{.2in}$s'$\hspace{.2in}$s$}
\vspace{0in}
\centerline{$\alpha$\hspace{.85in}$\alpha$\hspace{1.25in}}
\vspace{0in}
\centerline{\hspace{3.4in}$c'$}
\vspace{.7in}
\centerline{ $a'$\hspace{.75in}$w$\hspace{.75in}$b'$\hspace{1.25in}}
\vspace{.2in}
\caption{Unfolding surface.}
\label{unfold}
\end{figure}

Now, we have 
$$\d(b, b'c') - \d(b, a'b') = \d(b,b')(\sin\gamma-\sin\alpha).$$
By choosing $q$ large enough, the above difference can be made as large as 
wished; we make it larger than $\e$. Suppose now a (rectifiable)  path from
$x_0$ to $b'c'$ crosses $abc$ at $y$. Then its length is at least 
$$\d(x_0, y)+\d(y,b'c')\ge \delta(x_0, u)+\delta (b, b'c') \ge \d(x_0,u)+\d(b,a'b')+\e$$
$$>\d(x_0,u)+\d(b,a'b') \ge \d(x_0,u)+\d(u,a'b')\ge \d(x_0, a'b'),$$
so it cannot be a shortest path from $x_0$ to $a'b'c'$.
Analogously, no path from $x_0$ to $c'a'$ is a shortest path from $x_0$ to 
$a'b'c'$. This completes the proof of all desired properties for $P$. 
It is clear that, iterating this procedure, the shortest path from $x_0$ to 
the last constructed triangle has a steadily increasing spiralling number
(by a rate close to 1/3 for each new triangle). Thus
$s(P)$ can  be made as large as we wish.
The theorem is proved.
\end{proof}

\begin{remark} In the above construction the path $P$ can be forced to 
turn, each time, right or left,  as wished.
\end{remark}

In \cite{Pach}, J{\'a}nos Pach gives a neat example of a polytope
$K\subset {\mathbb R}^3$ and a shortest path $P$ on the boundary of $K$ with
non-smooth points $z_0,\ldots , z_n$ such that $\sum_{i=1}^n \gamma_i$ is not
bounded. Here $\gamma_i$ is the angle between the outer normals of the facets
containing the segments $[z_{i-1},z_i]$ and $[z_i,z_{i+1}]$. We mention that
our construction has the same property.

\end{document}